\input amstex
\magnification=\magstep1
\input epsf
\input amssym.def
\input amssym
\pageno=1
\baselineskip 14 pt
\def \pop#1{\vskip#1 \baselineskip}

\font\gr=cmbx12

\def \et{\operatorname {et}}

\def \Fr{\operatorname {Fr}}
\def \Spec{\operatorname {Spec}}
\def \Spf{\operatorname {Spf}}
\def \lim{\operatorname {lim}}

\def \fppf{\operatorname {fppf }}

\def \Spf{\operatorname {Spf}}

\def \Fr{\operatorname {Fr}}

\def \SS-Deg{\operatorname {SS-Deg}}
\def \SNS-Deg{\operatorname {SNS-Deg}}
 \def \DS-Deg{\operatorname {DS-Deg}}
\def \DNS-Deg{\operatorname {DNS-Deg}}
\def \F{\operatorname {F}}
\def \Id{\operatorname {Id}}

\par
\centerline
{\bf \gr On the existence of a torsor structure for Galois covers}

{\pop 2}
\noindent
\centerline {\bf \gr Mohamed Sa\"\i di}

{\pop 2}
\noindent

\centerline {\bf \gr Abstract}
{\pop 1}
Let $R$ be a complete discrete valuation ring with residue characteristic $p>0$. And let $f:Y\to X$ be a finite Galois cover
between flat and normal formal $R$-schemes of finite type with Galois group $G$ which is \'etale above the generic fibre of $X$ and 
such that the special fibre of $Y$ is reduced.
Assume that the special fibre of both $X$ and $Y$ are geometrically integral. If $G$ is cyclic of order $p$ then it is well known that
$f$ has the structure of a torsor under a finite and flat $R$-group scheme of rank $p$ (in the inequal characteristic case one has
to assume the existence of the $p$-th roots of unity in the ground field). In this note we give an example of a Galois cover 
as above with Galois group a cyclic group of order $p^2$ and such that $f$ doesn't have the structure of a  
torsor under a finite and flat $R$-group scheme of rank $p^2$ in the case where $R$ has equal characteristic $p$ (in the example one can even choose 
$X$ to be smooth). Such examples are certainly known to some 
experts but, to the best of my knowledge, do not appear in the literature. Also it is not surprising that such examples also exist in inequal characteristic.

\pop {1}
\par
\noindent
{\bf \gr I. Artin-Schreier-Witt theory of $p^n$-cyclic covers in 
characteristic $p$}.

\rm
\pop {.5}
\par
\noindent
{\bf \gr 1.1.}\rm\ In this section we review the Artin-Schreier-Witt theory (first developed in [W])
which provides explicit equations describing cyclic covers of degeree $p^n$ in 
characteristic $p$. We refeer the reader to the modern treatment of the theory as in [D-G].

\pop {.5}
\par
Let $X$ be a scheme of characteristic $p$ and denote by $X_{\et}$ the 
\'etale site on $X$. Let $n>0$ be an integer. We denote by $W_{n,X}$ 
(or simply $W_n$ if there is no confusion) the sheaf of Witt vectors of 
length $n$ on $X_{\et}$ (cf. [D-G], chapitre 5, 1). In the sequel any addition
or substarction of witt vectors will mean the addition and substarction
in the sens of Witt theory. We denote by $\F$ the 
Frobenius endomorphism of $W_{n,X}$ which is locally defined by
$F.(x_1,x_2,...,x_n)=(x_1^p,x_2^p,...,x_n^p)$, for a Witt vector
$(x_1,x_2,...,x_n)$ of length $n$, and by $\Id$ the identity 
automorphism of $W_{n,X}$. The following sequence is exact on $X_{\et}$:

$$(1)\ \ 0\to (\Bbb Z/p^n\Bbb Z)_X@>i_n>> W_{n,X} @> \F-\Id>> W_{n,X} \to 0$$
where $(\Bbb Z/p^n\Bbb Z)_X$ denotes the constant sheaf 
$(\Bbb Z/p^n\Bbb Z)$ on $X_{\et}$ and $i_n$ is the natural monomorphism 
which applies $1\in\Bbb Z/p^n\Bbb Z$ to $1\in W_n$ 
(cf [G-D], chapitre 5, 5.4). From the long cohomology exact sequence 
associated to (1) one deduces the following exact sequence:

$$(2)\ \ W_{n,X}(X)@>\F-\Id>>W_{n,X}(X)\to H^1_{\et}(X,\Bbb Z/p^n\Bbb Z)\to H^1_{\et}
(X,W_n)@>\F-\Id>> H^1_{\et}(X,W_n)$$

Assume that $X=\Spec A$ is affine in which case we have $H^1_{\et}(\Spec A,W_n)=0$ and hence 
an isomorphism:
$H^1_{\et}(\Spec A,\Bbb Z/p^n\Bbb Z)\simeq W_{n,A}(A)/(\F-\Id)(W_{n,A}(A))$. This 
isomorphism has the following interpretation: to an \'etale 
$\Bbb Z/p^n\Bbb Z$-torsor $f:Y\to X=\Spec A$ above $X$ 
corresponds a Witt vector $(a_1,a_2,...,a_n)\in W_{n,A}(A)$ 
of length $n$ which is uniquely determined 
modulo addition of elements of the form $F.(b_1,b_2,...,b_n)-(b_1,b_2,...,b_n)$. 
Moreover the equations 
$F.(x_1,x_2,...,x_n)-(x_1,x_2,...,x_n)=(a_1,a_2,...,a_n)$ where the $x_i$ 
are indeterminates are equations for the torsor $f$. More precisely there is 
a canonical factorisation of $f$ as $Y=Y_n@>f_{n}>> Y_{n-1}@>f_{n-1}>>...
@>f_2>> Y_1@>f_1>> Y_0:=X$ where each $Y_i=\Spec B_i$ is affine and $f_i:Y_i:=\Spec B_i\to f_{i-1}:
=\Spec B_{i-1}$ is the \'etale $\Bbb Z/p\Bbb Z$-torsor  corresponding to the algebra extension 
$B_{i+1}:=B_i[x_i]$. In the general case (where $X$ is not necessarily affine) the above equations provide local equations for an 
\'etale $\Bbb Z/p^n\Bbb Z$-torsor in characteristic $p$.

\pop {1}
\par
\noindent
{\bf \gr 1.2. Examples.}\rm\  In what follows $X$ is a scheme of characteristic $p$.

\pop {.5}
\par
\noindent
{\bf \gr 1.2.1. $\Bbb Z/p\Bbb Z$-Torsors.}\rm\ Let $f:Y\to X$ be an \'etale 
$\Bbb Z/p\Bbb Z$-torsor. Then locally $f$ is given by an equation $x^p-x=a$ 
where $a$ is a regular function on $X$ which is uniquely defined up to addition of 
elements of the form $b^p-b$ for some regular function $b$.

\pop {.5}
\par
\noindent
{\bf \gr 1.2.2. $\Bbb Z/p^2\Bbb Z$-Torsors.}\rm\ Let $f:Y\to X$ be an \'etale 
$\Bbb Z/p^2\Bbb Z$-torsor. Then we have a canonical factorisation of $f$ 
as: $Y=Y_2@>f_2>> Y_1@>f_1>>X$ where $f_2$ and $f_1$ are \'etale 
$\Bbb Z/p\Bbb Z$-torsors. The torsor $f$ is locally given by equations
of the form:
$$F.(x_1,x_2)-(x_1,x_2):=(x_1^p-x_1,x_2^p-x_2-p^{-1}\sum _{k=1}^{p-1}{p \choose k}x_1^{pk}(-x_1)^{p-k}) =(a_1,a_2)$$ 
for some regular functions $a_1$ and $a_2$ on $X$
and the Witt vector $(a_1,a_2)$ is uniquely determined up to addition (in the Witt
theory) of vectors of the form: 
$$(b_1^p,b_2^p)-(b_1,b_2):=(b_1^p-b_1,b_2^p-b_2-p^{-1}\sum _{k=1}^{p-1}{p \choose k}b_1^{pk}(-b_1)^{p-k})$$ 

Thus locally the torsor $f_1$ is defined by an equation: 
$$x_1^p-x_1=a_1$$ 
and $f_2$ by an equation: 
$$x_2^p-x_2=a_2+p^{-1}\sum _{k=1}^{p-1}{p \choose k}x_1^{pk}(-x_1)^{p-k}$$ 

Moreover if we replace the vector $(a_1,a_2)$ by
the vector $(a_1,a_2)+(b_1^p,b_2^p)-(b_1,b_2)$ the above equations are 
replaced by:  
$$x_1^p-x_1=a_1+b_1^p-b_1$$ 
and: 
$$x_2^p-x_2=a_2+b_2^p-b_2+p^{-1}\sum _{k=1}^{p-1}{p \choose k}x_1^{pk}(-x_1)^{p-k}-
p^{-1}\sum _{k=1}^{p-1}{p \choose k}b_1^{pk}(-b_1)^{p-k}-$$
$$p^{-1}\sum _{k=1}^{p-1}{p \choose k}(b_1^p-b_1)^{k}(a_1)^{p-k}$$ 
respectively.

\pop {.5}
\par
\noindent
{\bf \gr II. The group schemes $\Cal M_{n}$}\rm\ (cf. also [M], 3.2). 

\rm
\pop {.5}
\par
In all this paragraph we use the following notations: $R$ is a complete 
discrete valuation ring of equal characteristic $p$ with residue field $k$
and fraction field $K:=\Fr R$. We denote by $\pi$ a uniformising parameter 
of $R$.

\pop {.5}
\par
Let $n\ge 0$ be an integer and let
$\Bbb G_{a,R}=\Spec R[X]$ be the additive group scheme over $R$. The map: 
$$\phi_n:\Bbb G_{a,R}\to \Bbb G_{a,R}$$ 
given by: 
$$X\to X^p-\pi^{(p-1)n}X$$ 
is an isogeny of group schemes. The kernel of $\phi_n$ is denoted by $\Cal M_{n,R}:=\Cal M_{n}$. 
We have $\Cal M_{n}:=\Spec R[X]/(X^p-\pi^{(p-1)n}X)$, and
$\Cal M_{n}$ is a finite and flat $R$-group scheme of rank $p$. Further the following 
sequence is exact in the fppf topology:

$$(3)\ \ 0\to \Cal M_{n}\to \Bbb G_{a,R}@>\phi_n>>\Bbb G_{a,R}\to 0$$

If $n=0$ then the sequence $(3)$ is the Artin-Schreir sequence which is exact in the \'etale topology
and $\Cal M_{0}$ is the \'etale constant group scheme $(\Bbb Z/p\Bbb Z)_R$. If $n>0$ the sequence
$(3)$ has a generic fibre which is isomorphic to the \'etale Artin-Schreier sequence and a special fibre
isomorphic to the radicial exact sequence:

$$(4)\ \ 0\to \alpha_{p,k}\to \Bbb G_{a,k}@>x^p>>\Bbb G_{a,k}\to 0$$

Thus if $n>0$ the group scheme $\Cal M_{n}$ has a generic fibre which is \'etale isomorphic to
$(\Bbb Z/p\Bbb Z)_K$ and its special fibre is isomorphic to the infinitesimal group scheme 
$\alpha_{p,k}$. Let $X$ be an $R$-scheme. The sequence $(3)$ induces a long cohomology exact sequence:

$$(5)\ \ \Bbb G_{a,R}(X)@>\phi_n>>\Bbb G_{a,R}(X)\to H^1_{\fppf}(X,\Cal M_{n})
\to H^1_{\fppf}(X,\Bbb G_{a,R})@>\phi_n>> H^1_{\fppf}(X,\Bbb G_{a,R})$$

The cohomology group $H^1_{\fppf}(X,\Cal M_{n})$ classifies the isomorphism classes of
$\fppf$-torsors with group $\Cal M_{n}$ above $X$. The above sequence allows the 
following description of $\Cal M_{n}$-torsors:
locally a torsor $f:Y\to X$ under the group scheme $\Cal M_{n}$ is given by an equation 
$T^p-\pi^{(p-1)n}T=a$ where $T$ is an indeterminate and $a$ is a regular function on $X$ which is 
uniquely defined up to addition of elements of the form $b^p-\pi^{(p-1)n}b$ for
some regular function $b$. In particular if $H^1_{\fppf}(X,\Bbb G_{a,R})=0$ (e.g. if $X$ is affine) 
then an $\Cal M_{n}$-torsor above $X$ is globally defined by an equation as above.

\pop {.5}
\par
\noindent
{\bf \gr III. The example.}\rm \ 

\rm
\pop {.5}
\par
In all this paragraph we use the same notations as in II. Let $X:=\Spf A$ be a formal smooth affine $R$-scheme with geometrically connected fibres
(the assumption on $X$ being smooth is not essential for the example). Consider the cyclic $p^2$-cover $f:Y\to X$ given generically by the equation:
$$(T_1^p,T_2^p)-(T_1,T_2)=(\pi ^{-pm}a_1,a_2)$$
where $m=pm'$ is a positive integer divisible by $p$ and $a_1$ and $a_2$ are elements of $A$ such that the image 
$\bar a_1$ of $a_1$ modulo $\pi$ is not a $p$-power.
Then the generic fibre $f_K:Y_K\to X_K$ of $f$ is an \'etale $\Bbb Z/p^2\Bbb Z$ torsor. Further the finite cover $f$ factorises canonically
as $f_1\circ f_2$ where $f_2:Y\to Y_1$ and $f_1:Y_1\to Y$ are $p$-cyclic covers. The finite cover $f_1$ is given generically by the equation:
$$T_1^p-T_1=a_1\pi ^{-pm}$$
which can be transformed after making the change of variables $\tilde T_1:=\pi^{m}T_1$ to:
$$\tilde T_1^p-\pi ^{m(p-1)}\tilde T_1=a_1$$
which is a defining equation for the cover $f_1$. In particular we see that $f_1$ has the structure of a torsor under the $R$-group scheme $\Cal M_m$.
Its special fibre is the $\alpha_p$-torsor $f_{1,k}:Y_{1,k}\to X_{1,k}$ given by the equation $\tilde t_1^p=\bar a_1$ where $\tilde t_1$ is the image 
of $\tilde T_1$ modulo $\pi$. 

\pop {.5}
\par
The cover $f_2$ is generically given by the equation:
$$T_2^p-T_2=a_2+p^{-1}\sum _{k=1}^{p-1}{p \choose k}T_1^{pk}(-T_1)^{p-k}$$ 
and after adapting this to the change of variables $\tilde T_1:=\pi^{m}T_1$ we get the equation:
$$T_2^p-T_2=a_2+p^{-1}\sum _{k=1}^{p-1}{p \choose k}\tilde T_1^{pk}(-\tilde T_1)^{p-k}\pi^{-m(pk+p-k)}$$
which can be transformed after the change of variables $\tilde T_2:=\pi ^{\tilde m}T_2$ where $\tilde m:=m'(p(p-1)+1)$ to:
$$\tilde T_2^p-\pi ^{\tilde m(p-1)}\tilde T_2=\pi ^{m(p(p-1)+1)}a_2+p^{-1}\sum _{k=1}^{p-1}{p \choose k}\tilde T_1^{pk}(-\tilde T_1)^{p-k}
\pi^{\tilde m p-m(pk+p-k)}$$

The above equation is an integral equation defining the cover $f_2$. In particular one sees that $f_2$ has the structure of a torsor under the $R$-group 
scheme $\Cal M_{\tilde m}$. Its special fibre is the $\alpha_p$-torsor $f_{2,k}:Y_{2,k}\to X_{2,k}$ given by the equation 
$\tilde t_2^p=-\tilde t_1^{p(p-1)+1}$ where $\tilde t_2$ is the image of modulo $\pi$. In particular the special fibre $Y_k$ of $Y$ is reduced. 
Although both $f_1$ and $f_2$ have the structure of a torsor the next proposition shows that the cover $f$ doesn't.

\pop {.5}
\par
\noindent
{\bf \gr 3.1. Proposition.}\rm\ {\sl With the same notations as above and assume that $p>2$. Then the finite cover $f$
doesn't have the structure of a torsor under a finite and flat $R$-group scheme of rank $p^2$.}

\pop {.5}
\par
\noindent
{\bf \gr Proof.}\rm\ Let $b_1\in A$. Then the equations:
$(S_1^p,S_2^p)-(S_1,S_2)=(\pi^{-pm}a_1,a_2)+(\pi^{-pm}b_1^p,0)-(\pi^{-m}b_1,0)$ are also defining equations for the torsor $f_K$. 
The cover $f$ is then given by the equations (1):
$$\Tilde T_1^p-\pi ^{m(p-1)}\Tilde T_1=a_1$$ 
and  
$$\Tilde T_2^p-\pi ^{\tilde m(p-1)}\Tilde T_2=
\pi^{\tilde mp}a_2+p^{-1}\sum _{k=1}^{p-1}{p \choose k}\Tilde T_1^{pk}(-\Tilde T_1)^{p-k}\pi ^{\tilde mp-pk-p+k}$$ 

And $f$ is also given by the equations (2): 
$$\Tilde S_1^p-\pi ^{m(p-1)}\Tilde S_1=a_1+b_1^p-\pi^{m(p-1)}b_1$$ 
and: 
\newline
$S_2^p-S_2=a_2+p^{-1}(\sum _{k=1}^{p-1}{p \choose k}\Tilde S_1^{pk}(-\Tilde S_1)^{p-k}\pi ^{-m(pk+p-k)}-
\sum _{k=1}^{p-1}{p \choose k}b_1^{pk}(-b_1)^{p-k}\pi ^{-m(pk+p-k)}$
\newline
$-\sum _{k=1}^{p-1}{p \choose k}
(b_1^p-\pi^{m(p-1)}b_1)^{k}(a_1)^{p-k}\pi^{-mp^2})$. 
By reducing modulo $\pi$ both equations (1) and (2) we obtain equations (1)' and (2)' for the cover $f_k$ on the level of special fibres. 
Suppose that $f$ has the structure of a 
torsor under a finite and flat $R$-group scheme $G_R$ of rank $p^2$. Then the special fibre $f_k$ of $f$ is a torsor under the special fibre $G_k$ of $G$ 
which is a group scheme of rank $p^2$ and an extension of $\alpha_p$ by $\alpha_p$. In particular both equations (1)' and (2)' that we obtain for $f_k$ must 
define the same $\alpha_p$-torsor $f_{1,k}$ and $f_{2,k}$. This is the case for $f_{1,k}$ since both equations are $\tilde t_1^p=\bar a_1$ and 
$\tilde s_1^p=\bar a_1+\bar b_1^p$ but we will see that this is not the case 
for $f_{2,k}$ if $p>2$. Indeed one can see after some (easy) computations that the equations we obtain for $f_{2,k}$ are: 
$\tilde t_2^p=-\tilde t_1^{p(p-1)+1}$ and $\tilde s_2^p=-\tilde s_1^{p(p-1)+1}+\bar b_1^{p(p-1)+1}+\sum _{k=1}^{p-1}(\bar b_1^{p(k-1)+1}(\tilde s_1-\bar 
b_1)^{p(p-k)}+\bar b_1^{pk}(\tilde s_1-\bar b_1)^{p(p-k-1)+1}))$, where $\tilde s_2$ is the image of $\pi^{\tilde m}S_2$ modulo $\pi$,
which gives $\tilde s_2=-\tilde s_1^{p(p-1)+1}+\bar b_1^{p(p-1)+1}+((\tilde s_1-\bar 
b_1)^{p-1}+\bar b_1^{p-1})\sum _{k=1}^{p-1}\bar b_1^{p(k-1)+1}(\tilde s_1-\bar b_1)^{p(p-k-1)+1}$. An easy verification shows that these two equations define 
non isomorphic $\alpha_p$-torsors (although they define isomorphic covers). 
If for example $p=3$ then the two equations are: $\tilde t_2^3=-\tilde t_1^7$ and $\tilde s_2^3=-\tilde s_1^7+\bar b_1^7+
((\tilde s_1-\bar b_1)^2-\bar b_1^2)(\bar b_1(\tilde s_1-\bar b_1)^4+\bar b_1^4(\tilde s_1-\bar b_1))$. Using that $\tilde s_1=\tilde t_1+\bar b_1$ we get 
$\tilde s_2^3= -(\tilde t_1+\bar b_1)^7+\bar b_1^7+(\tilde t_1^2-\bar b_1^2)(\bar b_1\tilde t_1^4+\bar b_1^4\tilde t_1)$ thus 
$\tilde s_2^3=-\tilde t_1^7-34\tilde t_1^3\bar b_1^4-8\tilde t_1\bar b_1^6=-\tilde t_1^7+
2\tilde t_1^3\bar b_1^4+\tilde t_1\bar b_1^6$ and since $2\tilde t_1^3\bar b_1^4+\tilde t_1\bar b_1^6=2\bar a_1\bar b_1^4+\tilde t_1\bar b_1^6$ 
is not a cube we deduce that the two equations do not define the same $\alpha_3$-torsor.

\pop {2}
\par
\noindent
{\bf \gr References}
\rm

\pop {.5}
\par
\noindent
[D-G] M. Demazure and P. Gabriel, Groupes Alg\'ebriques, Tome 1, Masson and CIE \'Editeur, Paris, North-Holland Publishing Company, 
Amsterdam, (1970).

\pop {.5}
\par
\noindent
[M] S. Maugeais, {\sl Rel\`evement des rev\^etements $p$-cycliques des courbes rationnelles semi-stable},
Math. Ann., 327, 365-393, (2003).

\pop {.5}
\par
\noindent
[W] E. Witt, {\sl Zyclische K\"orper und algebren der characteristic $p$ vom 
grad $p^n$. Struktur diskret bewerteter perfekter K\"orper mit vollkommenem 
Restklassenk\"orper der Characteristic $p$}, J. F\"ur die reine und angewandte
Mathematik (Crelle), 176, 126-140, (1937).

\enddocument